\documentclass{article}
\usepackage{epsfig,subfigure,extpfeil}
\usepackage{amsmath,amsthm}
\usepackage{amssymb}
\usepackage{graphicx,epstopdf,latexsym,amsfonts,fancyhdr}
\usepackage{booktabs,arydshln,multirow}
\usepackage{setspace}
\usepackage{xcolor}
\usepackage{stmaryrd}
\usepackage{cite}
\usepackage[utf8]{inputenc}
\newtheorem{theorem}{Theorem}[section]

\newtheorem{lemma}[theorem]{Lemma}

\newtheorem{remark}[theorem]{Remark}

\numberwithin{equation}{section} \topmargin=-2.5cm \oddsidemargin=0.4cm
\newcommand{\normmm}[1]{{\left\vert\kern-0.25ex\left\vert\kern-0.25ex\left\vert #1
    \right\vert\kern-0.25ex\right\vert\kern-0.25ex\right\vert}}
\begin{document}
\title{Construction of high-order robust theta-methods with applications in anomalous models
\thanks{Corresponding author.
\newline \emph{Email addresses:} baolimath@126.com
}}
\date{ }
\author{Baoli Yin${}^{1*}$, Guoyu Zhang${}^1$, Yang Liu${}^1$, Hong Li${}^1$
\\\small{\emph{${}^1$School of Mathematical Sciences, Inner Mongolia University, Hohhot 010021,
China;}}
}
\date{}
 \maketitle
  {\color{black}\noindent\rule[0.5\baselineskip]{\textwidth}{0.5pt} }
\noindent \textbf{Abstract:}
A general conversion strategy by involving a shifted parameter $\theta$ is proposed to construct high-order accuracy difference formulas for fractional calculus operators.
By converting the second-order backward difference formula with such strategy, a novel $\theta$-scheme with correction terms is developed for the subdiffusion problem with nonsmooth data, which is robust even for very small $\alpha$ and can resolve the initial singularity.
The optimal error estimates are carried out with essential arguments and are verified by numerical tests.
\\
 {\color{black}\noindent\rule[0.5\baselineskip]{\textwidth}{0.5pt} }
\def\REF#1{\par\hangindent\parindent\indent\llap{#1\enspace}\ignorespaces}
\newcommand{\h}{\hspace{1.cm}}
\newcommand{\hh}{\hspace{2.cm}}
\newtheorem{yl}{\hspace{0.cm}Lemma}
\newtheorem{dl}{\hspace{0.cm}Theorem}
\newtheorem{re}{\hspace{0.cm}Remark}
\renewcommand{\sec}{\section*}
\renewcommand{\l}{\langle}
\renewcommand{\r}{\rangle}
\newcommand{\be}{\begin{eqnarray}}
\newcommand{\ee}{\end{eqnarray}}
\normalsize \vskip 0.2in
\section{Introduction}\label{sec.intr}
The subdiffusion transport mechanism in recent years has received much attention for the fact that some physical processes including the electron transport, thermal diffusion, and protein transport, among others, reveal that the underlying stochastic process is the continuous time random walk instead of the Brownian motion \cite{JinLazarovZhouconcise,MetzlerKlafterrandom}.
In this study, we develop robust time-stepping methods for the following $\alpha$th ($\alpha \in (0,1)$) order subdiffusion problem
\begin{equation}\label{I.1}\begin{split}
\begin{cases}
  \partial_t^\alpha u(\boldsymbol x,t)-\Delta u(\boldsymbol x,t)=f(\boldsymbol x,t), & (\boldsymbol x,t)\in \Omega\times (0,T], \\
  u(\boldsymbol x,t)=0, & \boldsymbol x \in \partial\Omega, ~t\in (0,T],\\
  u(\boldsymbol x,0)=v(\boldsymbol x), & \boldsymbol x \in \Omega,
\end{cases}
\end{split}\end{equation}
where the space $\Omega \in \mathbb{R}^d$ $(d=1,2,3)$ is a bounded convex polygonal domain with the boundary denoted by $\partial \Omega$.
The operator $\Delta: D(\Delta) \to L^2(\Omega)$ stands for the Laplacian with $D(\Delta)=H_0^1(\Omega)\cap H^2(\Omega)$, and $f:(0,T]\to L^2(\Omega)$ is a given function.
The initial function $v$, depending on its smoothness, belongs to $D(\Delta)$ or $L^2(\Omega)$.
$\partial_t^\alpha$ is the Caputo fractional operator satisfying $\partial_t^\alpha \phi=D_t^\alpha (\phi-\phi(0))$ for $\alpha \in (0,1)$, where $D_t^\alpha$, known as the Riemann-Liouville fractional operator, is defined by
\[
(D_t^\alpha \phi)(t)=\frac{1}{\Gamma(1-\alpha)}\frac{\mathrm{d}}{\mathrm{d}t}
\int_{0}^{t}\frac{\phi(s)}{(t-s)^\alpha}\mathrm{d}s.
\]
\par
The literature on subdiffusion is vast, for example, the solution regularity exploration can be found in \cite{SakamotoYamamoto}, and some numerical studies were developed in \cite{YinLiuLiZhang1,JinLiZhou1,JinLiZh2,WangWangYin,GaoSun,LiZhaoChen,LiaoMcLeanZhang}, to mention just a few.
See also the overview article \cite{JinLazarovZhouconcise}.
It is well known that the problem (\ref{I.1}) is characterized by the initial singularity of its solution, which frustrates most high-order numerical methods in case the singularity is overlooked.
In \cite{YinLiuLiZhang1}, we proposed a modified $\theta$-method which can preserve the optimal accuracy for $\theta \in (0,\frac{1}{2})$.
As mentioned in \cite{YinLiuLiZhang1}, the case $\theta=\frac{1}{2}$ has deserved much more our attention since the correction terms vanish when $\theta=\frac{1}{2}$, enlightening us that a carefully designed time-stepping method should automatically resolve the singularity.
To sum up, our contribution in this study is twofold:
\begin{itemize}
  \item A novel strategy is developed which can transfer known time-stepping methods such as the fractional BDF2 to more robust methods.
  \item Rigorous arguments of the optimal error estimates of the transformed fractional BDF2 are provided for the subdiffusion problem (\ref{I.1}).
\end{itemize}
\par
The rest of the article is outlined as follows.
In section \ref{sec.Novel}, a novel strategy is proposed to introduce a shifted parameter $\theta$ into known stepping methods, based on which the fully discrete scheme for (\ref{I.1}) is constructed.
In section \ref{sec.Optimal}, the rigorous error estimates are provided and their correctness is fully validated in section \ref{sec.tests}.
Finally, some concluding remarks are made in section \ref{sec.conc}.
\section{Novel $\theta$-schemes}\label{sec.Novel}
We first propose some general results on constructing high-order accuracy difference formulas for fractional calculus based on generating function (GF) reformulation.
Assume $\varpi_p(\zeta)$ is a GF of the convolution quadrature (CQ) \cite{Lubich} with convergence order $p$, and let $\delta(\zeta)=\sum_{j=1}^{p}\frac{1}{j}(1-\zeta)^j$ denote the GF of backward difference formulas (BDF) with $p\leq 6$.
\begin{lemma}\label{lem.novel.1} (General conversion strategy)
Define $\omega(\zeta)=\varpi_p(\zeta) e^{\theta \delta(\zeta)}$, $\theta \in \mathbb{R}$,
then $\omega(\zeta)$ can generate a $\theta$-method which is convergent of order $p$.
\begin{proof}
The function $e^{\theta \delta(\zeta)}$ is sufficiently differentiable on the unit circle and thus its Fourier coefficients decay faster than, e.g., $O(n^{-k})$ for any positive integer $k$.
Then the asymptotic property of $\omega_n$ is fully determined by $\varpi_n$ which, by the stability in CQ (i.e., $\varpi_n=O(n^{-\alpha-1})$, see Definition 2.1 in \cite{Lubich}), leads to $\omega_n=O(n^{-\alpha-1})$.
Moreover, by the consistency of $\varpi_p(\zeta)$ (see Definition 2.2 in \cite{Lubich})) and the backward difference formulas, i.e.,
\[
\tau^{-\alpha}\varpi_p(e^{-\tau})=1+O(\tau^p),\quad
\tau^{-1}\delta(e^{-\tau})=1+O(\tau^p),
\]
we have
$
\tau^{-\alpha}e^{\theta \tau}\omega(e^{-\tau})
=\tau^{-\alpha}\varpi_p(e^{-\tau})e^{\theta \tau}e^{\theta \delta(e^{-\tau})}=1+O(\tau^p),
$
indicating that $\omega(\zeta)$ is consistent of order $p$ which, combined with $\omega_n=O(n^{-\alpha-1})$, completes the proof of the lemma (see Theorem 1 in \cite{LiuYinunified}).
\end{proof}
\end{lemma}
\begin{remark}\label{rem.1}
Lemma \ref{lem.novel.1} indicates we can approximate $\phi(t_{n-\theta})$ by a discrete convolution as
\begin{equation}\label{novel.1.2}
\sum_{j=0}^{n}\theta_j\phi(t_{n-j}),\quad \text{where $\theta_j$ is generated by~} \sum_{j=0}^{\infty}\theta_j \zeta^j=e^{\theta \delta(\zeta)}.
\end{equation}
\end{remark}
\begin{lemma}\label{lem.novel.2}
Assume $\omega(\zeta)$ takes the form $\big[P(\zeta)\big]^\alpha e^{\theta Q(\zeta)}$ where $P(\zeta)$ and $Q(\zeta)$ are polynomials such that $\omega(\zeta)$ is analytic within the open unit disc, then
\begin{equation}\label{novel.2}
\omega_n=\frac{1}{nP(0)}\bigg[\omega_0 G_{n-1}+\sum_{k=1}^{n-1}\omega_{n-k}\big(G_{k-1}-(n-k)P_k\big)\bigg], \quad n \geq 1,\quad
\omega_0=\big[P(0)\big]^\alpha e^{\theta Q(0)},
\end{equation}
where $G_k$ is the coefficients of $G(\zeta)$ defined by $G(\zeta)=\alpha P'(\zeta)+\theta P(\zeta)Q'(\zeta)$.
\begin{proof}
Take the derivative of $\omega(\zeta)=\big[P(\zeta)\big]^\alpha e^{\theta Q(\zeta)}$ w.r.t $\zeta$ and multiply both sides by $P(\zeta)$ to obtain
\[
P(\zeta)\omega'(\zeta)=\omega(\zeta)G(\zeta).
\]
The formula (\ref{novel.2}) then follows by taking the $n$th coefficient of both sides of the above equality.
\end{proof}
\end{lemma}
\par
It is notable that the algorithm (\ref{novel.2}) is efficient since $G(\zeta)$ and $P(\zeta)$ have finitely many nonzero coefficients, and thus the computing complexity to obtain $\{\omega_j\}_{j=0}^N$ is of $O(N)$.
\par
Denote by $u^n$ the approximation to $u(t_n)$, and introduce the symbols for general functions $\phi$
\begin{equation}\label{novel.2.0}
\phi^{n-\theta}=\sum_{j=0}^{n}\theta_j \phi^{n-j},\quad
D_\tau^{\alpha,n-\theta}\phi=\tau^{-\alpha}\sum_{j=0}^{n}\omega_j \phi^{n-j}
\end{equation}
where $\theta_j$ is defined in (\ref{novel.1.2}) with $\delta(\zeta)=\frac{3}{2}-2\zeta+\frac{1}{2}\zeta^2$, and $\omega_j$ is generated by $\omega(\zeta)=\big[\delta(\zeta)\big]^\alpha e^{\theta\delta(\zeta)}$.
In accordance with Lemma \ref{novel.1.2} (see also Remark \ref{rem.1}), $\phi^{n-\theta}$ and $D_\tau^{\alpha,n-\theta}\phi$ both are of second-order accuracy to their continuous counterparts.
To formulate the fully discrete scheme of the model, define the finite element space as
$
V_h=\{\chi_h \in H_0^1(\Omega): \chi_h |_e \text{ is a linear polynomial function},~ e \in \mathcal{T}_h\}
$
where $\mathcal{T}_h$ is a shape regular, quasi-uniform triangulation of $\Omega$.
\par
Let $P_h:L^2(\Omega)\to V_h$ and $R_h:H_0^1(\Omega)\to V_h$ stand for the $L^2(\Omega)$ and Ritz projection, respectively, and define $\Delta_h: V_h\to V_h$ as the discrete Laplacian.
By replacing $u(t)$ with $w(t)+v$ and $f(t)$ with $g(t)+f(0)$ in (\ref{I.1}), the space semi-discrete scheme then reads
\begin{equation}\label{novel.2.1}
D_t^\alpha w_h(t)-\Delta_h w(t)=g_h(t)+f_h^0+\Delta_h v_h,
\end{equation}
where $g_h:=P_h g$, $f_h^0=P_h f(0)$ and $v_h=R_h v$ if $v\in D(\Delta)$ or $v_h=P_h v$ if $v \in L^2(\Omega)$.
Then the fully discrete scheme can be stated as finding $W_h^n\in V_h$ such that
\begin{equation}\label{novel.3}
D_\tau^{\alpha,n-\theta}W_h-\Delta_h W_h^{n-\theta}=g_h^{n-\theta}+f_h^0+\Delta_h v_h,\quad n\geq 1, \quad \theta \in (-1,1).
\end{equation}

\par
In general cases, the scheme (\ref{novel.3}) can only result first-order convergence rate at positive time due to the initial singularity of the solution.
We propose a corrected scheme, with the motivation explained in the next section, by resorting to a single-step modification:
\begin{equation}\label{novel.4}\begin{split}
D_\tau^{\alpha,1-\theta}W_h-\Delta_h W_h^{1-\theta}=(\theta+3/2)(\Delta_h v_h+f_h^0)
+g_h^{1-\theta},\quad n=1,
\\
D_\tau^{\alpha,n-\theta}W_h-\Delta_h W_h^{n-\theta}=g_h^{n-\theta}+f_h^0+\Delta_h v_h,\quad n\geq 2.
\end{split}\end{equation}
\par
We note that for $\theta=-\frac{1}{2}$, the scheme (\ref{novel.4}) recovers exactly (\ref{novel.3}), indicating that (\ref{novel.3}) can resolve the initial singularity automatically if the problem is discretized at the point $t_{n+\frac{1}{2}}$.
\section{Optimal error estimates}\label{sec.Optimal}
The error estimate is based on solution representation and estimates of some kernels.
Denote by $\widehat{\phi}$ the Laplace transform of $\phi$.
Then, using the Laplace transform and its inverse transform, we obtain
\begin{equation}\label{Optimal.1}\begin{split}
w_h(t)&=-\frac{1}{2\pi{\rm i}}\int_{\Gamma_{\sigma,\epsilon}}e^{zt}\big[K(z)(\Delta_h v_h+f_h(0))+zK(z)\widehat{g_h}(z)\big]\mathrm{d}z,
\end{split}\end{equation}
where $K(z)=-z^{-1}(z^\alpha-\Delta_h)^{-1}$ stands for the kernel function, and the contour (with the direction of an increasing imaginary part) $\Gamma_{\sigma,\epsilon}$ is defined by
\[
\Gamma_{\sigma,\epsilon}:=\{z\in\mathbb{C}:|z|=\epsilon, |\arg z|\leq \sigma\}\cup\{z\in\mathbb{C}: z=re^{\pm{\rm i}\sigma}, r\geq \epsilon\}.
\]
\begin{theorem}\label{thm.1}
For $\alpha \in (0,1)$ and $\theta \in (-1,1)$, there exist $\sigma_0 \in (\pi/2,\pi)$ and $\epsilon_0>0$ both of which are free of $\alpha$ and $\tau$ such that for any $\sigma\in (\pi/2,\sigma_0)$ and any $\epsilon<\epsilon_0$, the solution of (\ref{novel.4}) takes the form
\begin{equation}\label{Optimal.2}\begin{split}
W_h^n=-\frac{1}{2\pi{\rm i}}\int_{\Gamma^\tau_{\sigma,\epsilon}}e^{zt_n}\big[\ell(e^{-z\tau})K(\delta_\tau(e^{-z\tau}))
(\Delta_h v_h+f_h^0)
+\tau\delta_\tau(e^{-z\tau})K(\delta_\tau(e^{-z\tau}))g_h(e^{-z\tau})\big]\mathrm{d}z,
\end{split}\end{equation}
where $\Gamma^\tau_{\sigma,\epsilon}=\{z\in \Gamma_{\sigma,\epsilon}:|\Im(z)|\leq \pi/\tau\}$,
$\delta_\tau(\zeta)=\delta(\zeta)/\tau$ and $\ell(\zeta)=\delta(\zeta)\zeta\big(\frac{1}{1-\zeta}+\theta+\frac{1}{2}\big)e^{-\theta\delta(\zeta)}$.
\begin{proof}
Multiply both sides of (\ref{novel.4}) by $\zeta^n$ and sum the index $n$ from $1$ to $\infty$ to yield
\[
\sum_{n=1}^{\infty}\zeta^n D_\tau^{\alpha,n-\theta}W_h
-\sum_{n=1}^{\infty}\zeta^n\Delta_h W_h^{n-\theta}
=\sum_{n=1}^{\infty}\zeta^n g_h^{n-\theta}
+(f_h^0+\Delta_h v_h)\bigg(\sum_{n=1}^{\infty}\zeta^n+(\theta+1/2)\zeta\bigg),
\]
which, by definitions of symbols in (\ref{novel.2.0}), leads to
\[
\big(\big[\delta_\tau(\zeta)\big]^\alpha
-\Delta_h\big) W_h(\zeta)
=g_h(\zeta)
+(f_h^0+\Delta_h v_h)\kappa(\zeta),
\]
where $\kappa(\zeta)=\zeta\big(\frac{1}{1-\zeta}+\theta+\frac{1}{2}\big)e^{-\theta\delta(\zeta)}$.
By Lemma B.1 in \cite{JinLiZhou1}, for fixed constant $\phi_0\in (\pi/2,\pi)$, there exists $\sigma_0\in (\pi/2,\pi)$ which depends only on $\phi_0$, for any $\sigma\in (\pi/2,\sigma_0)$ and any $\epsilon<\epsilon_0$ where $\epsilon_0$ is small enough, $\delta_\tau(e^{-z\tau})|_{z\in \Gamma_{\sigma,\epsilon}^\tau}\in \Sigma_{\phi_0}:=\{z\in \mathbb{C}:|\arg z|<\phi_0, z\neq 0\}$.
By Cauchy integral formula, we have the expression for $W_h^n$ by
\[
W_h^n=\frac{1}{2\pi{\rm i}}\int_{|\zeta|=\varepsilon}\frac{W_h(\zeta)}{\zeta^{n+1}}\mathrm{d}\zeta
\xlongequal{\zeta=e^{-z\tau}}\frac{\tau}{2\pi{\rm i}}\int_{\Gamma_\varepsilon^\tau}e^{zt_n}W_h(e^{-z\tau})\mathrm{d}z
\]
where $\Gamma^\tau_\varepsilon:=\big\{z=-\frac{1}{\tau}\ln\varepsilon+{\rm i}y: y\in \mathbb{R}, |y|\leq \pi/\tau\big\}$.
Let $\mathcal{L}$ be the region enclosed by contours $\Gamma_{\sigma,\epsilon}^\tau$, $\Gamma_\varepsilon^\tau$, $\Gamma^\tau_{\pm}:=\mathbb{R}\pm {\rm i}\pi/\tau$ (oriented from left to right), one can check $W_h(e^{-z\tau})$ is analytic for $z\in \overline{\mathcal{L}}$.
By using the Cauchy integral formula again, and noting that the integral values along $\Gamma_-^\tau$ and $\Gamma_+^\tau$ are opposite, the result (\ref{Optimal.2}) follows readily by taking $\ell(\zeta)=\tau\delta_\tau(\zeta)\kappa(\zeta)$.
The proof is completed.
\end{proof}
\end{theorem}
\begin{remark}\label{rem.2}
The arguments for Theorem \ref{thm.1} reveal the superiority of our scheme that, on the one hand for arbitrary $\theta$, the transform function $e^{-\theta \delta(\zeta)}|_{\zeta=e^{-z\tau}}$ appeared in $\kappa(\zeta)$ is analytic for $z\in \overline{\mathcal{L}}$, in contrast to the transform function $\frac{1}{1-\theta+\theta\zeta}|_{\zeta=e^{-z\tau}}$ in \cite{YinLiuLiZhang1} which is singular at points $z=\pm \frac{\pi}{\tau}\in \overline{\mathcal{L}}$ when $\theta=\frac{1}{2}$ (in which case, the Crank-Nicolson scheme is excluded).
See also \cite{JinLiZh2,WangWangYin} for similar situations.
Therefore, our scheme or numerical analysis is robust against the shifted parameter $\theta$.
On the other hand, thanks to Lemma \ref{lem.novel.1}, the function $\delta_\tau(\zeta)$ appeared in (\ref{Optimal.2}) is independent of $\alpha$, allowing us to develop robust analysis even for small $\alpha$.
We argue that such kind of robustness is not available for schemes in \cite{JinLiZh2,WangWangYin,YinLiuLiZhang1} as $\delta_\tau(\zeta)$ in those schemes are singular at $\alpha=0$, leading to the blow-up of constants $C$ in their estimates.
See \textit{Example 2} in section \ref{sec.tests}.
\end{remark}
\begin{lemma}\label{lem.optimal.1}
Let $\Gamma^\tau_{\sigma,\epsilon}$ be the contour defined in Theorem \ref{thm.1}.
For given $\theta \in (-1,1)$ and any $z\in \Gamma^\tau_{\sigma,\epsilon}$, there holds
\begin{equation}\label{Optimal.2.1}
|\ell(e^{-z\tau})-1|\leq C\tau^2 |z|^2,
\end{equation}
where $C$ is independent of $\tau, z$, but may dependent on $\theta$.
\begin{proof}
Since $|z|\tau \leq \pi/\sin\sigma<+\infty$, we only need to prove (\ref{Optimal.2.1}) for sufficiently small $|z|\tau$.
By the expansion of $\ell(\zeta)$ at the point $\zeta=1$, we have
$
\ell(\zeta)=1+c(\theta)(1-\zeta)^2+(1-\zeta)^3 r(\zeta),
$
where $r(\zeta)$ is analytic at $\zeta=1$.
One then immediately gets
$
\ell(e^{-z\tau})=1+c(\theta)\tau^2|z|^2+o(\tau^2|z|^2)
$,
which completes the proof of the lemma.
\end{proof}
\end{lemma}

\begin{theorem}\label{thm.2}
Suppose $u_h(t):=w_h(t)+v_h$ is the solution of the space semi-discrete scheme of (\ref{I.1}), and $U_h^n:=W_h^n+v_h$ is the solution of the fully discrete scheme of (\ref{I.1}).
If $f\in W^{1,\infty}(0,T;L^2(\Omega))$ and $\int_{0}^{t}(t-s)^{\alpha-1}\|f''(x)\|\mathrm{d}s\in L^{\infty}(0,T)$ where $\|\cdot\|$ denotes the $L^2$ norm, then
\begin{equation}\label{Optimal.3}
\|U_h^n-u_h(t_n)\|\leq C\tau^2\bigg(\mathcal{R}(t_n,v)
+t_n^{\alpha-2}\|f(0)\|
+t_n^{\alpha-1}\|f'(0)\|
+\int_{0}^{t_n}(t_n-s)^{\alpha-1}\|f^{''}(s)\|\mathrm{d}s
\bigg),
\end{equation}
where $\mathcal{R}(t_n,v)=t_n^{\alpha-2}\|\Delta v\|$ if $v \in D(\Delta)$ and $\mathcal{R}(t_n,v)=t_n^{-2}\|v\|$ if $v \in L^2(\Omega)$.
The constant $C$ is independent of $\tau,\alpha,n,N$ and $f$, but may depend on $\theta$.
\begin{proof}
The arguments for this theorem is essentially based on Lemma \ref{lem.optimal.1} and the following estimates on $\delta_\tau(\zeta)$, which can be found in \cite{JinLiZhou1},
\[
|\delta_\tau(e^{-z\tau})-z|\leq C\tau^2|z|^3,\quad
|\delta_\tau^\alpha(e^{-z\tau})-z^\alpha|\leq C\tau^2|z|^{2+\alpha},\quad
C_1|z|\leq |\delta_\tau(e^{-z\tau})|\leq C_2|z|.
\]
Then, the result (\ref{Optimal.3}) is followed after a lengthy but standard analysis for the contour integral, which is omitted for space reasons.
\end{proof}
\end{theorem}
\begin{remark}
The error $u-u_h$ of the space semi-discrete scheme (\ref{novel.2.1}) has been well studied by researchers which is not our main concern in this article.
Interested readers can refer, e.g., \cite{JinLazarovZhou} for more information.
\end{remark}
\section{Numerical tests}\label{sec.tests}
\textit{Example 1.} Let $T=1$.
Depending on the smoothness of $v$, we consider two cases:
\par
~(i) $f=0$, $v=\sin x \in D(\Delta)$, $\Omega=(0,\pi)$, with the exact solution $u(x,t)=E_\alpha(-t^\alpha)\sin x$;
\par
(ii) $f=0$, $v=\chi_{(0,1/2)}$, $\Omega=(0,1)$;
\par
\begin{figure}[htbp]
\centering

\subfigure[]{
\begin{minipage}[t]{0.48\linewidth}
\centering
\includegraphics[width=1\textwidth]{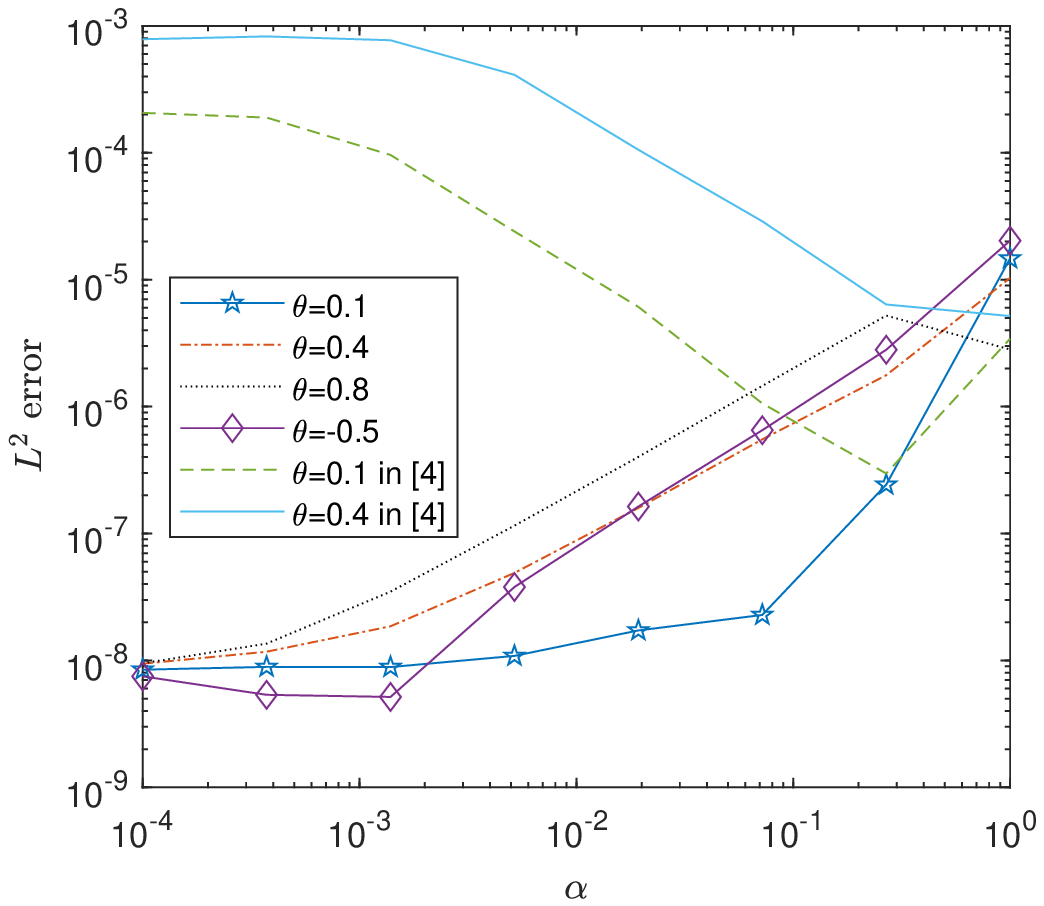}
\end{minipage}%
}%
\subfigure[]{
\begin{minipage}[t]{0.48\linewidth}
\centering
\includegraphics[width=1\textwidth]{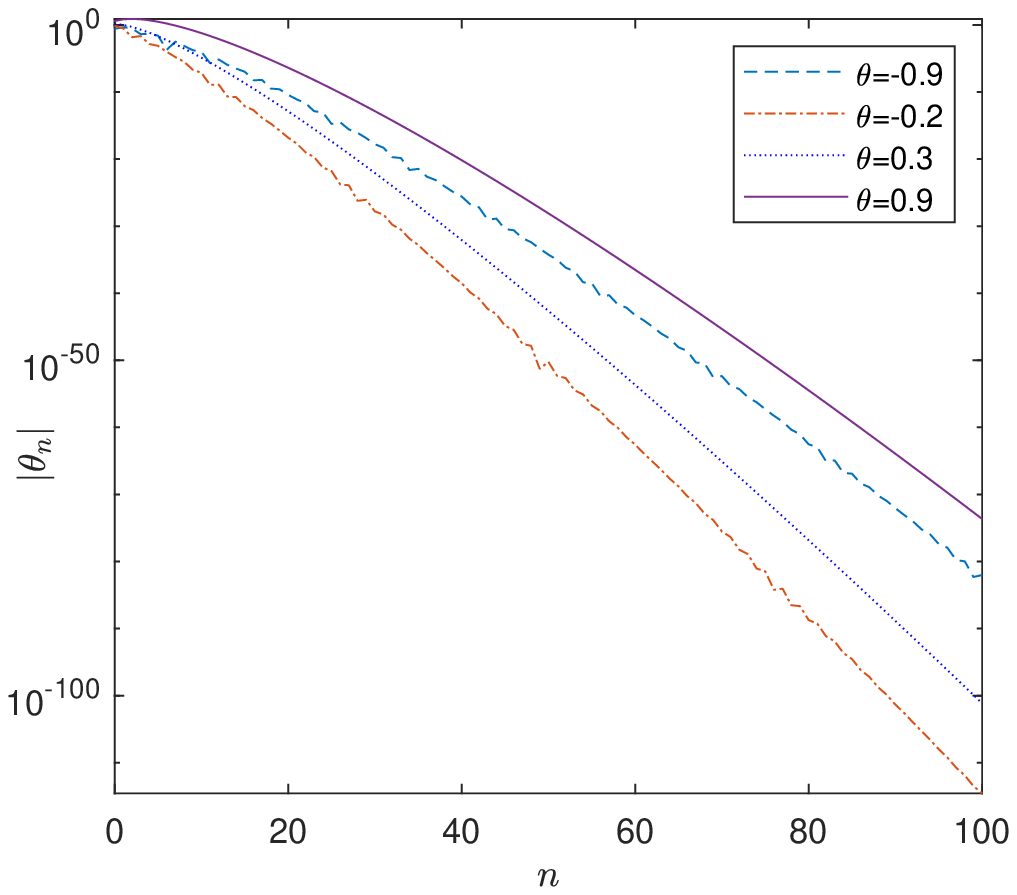}
\end{minipage}%
}%
\centering
\caption{(a) Comparison of $L^2$ error between our scheme and that in \cite{YinLiuLiZhang1} for different $\alpha$. (b) Exponential decay of the weights $|\theta_n|$ defined in (\ref{novel.1.2}).}\label{C1}
\end{figure}

In Table \ref{tab1} and Table \ref{tab2}, we present the $L^2$ error and convergence rates for different $\alpha$ and $\theta$ for schemes (\ref{novel.3}) and (\ref{novel.4}), respectively.
One observes that the scheme (\ref{novel.4}) with correction terms results in optimal convergence rates while the scheme (\ref{novel.3}) is of first-order accuracy except for $\theta=-0.5$, both of which are in line with our theoretical results.
\begin{table}[htbp]
\centering
\caption{$L^2$ error and convergence rates at time $t=0.5$ of \textit{Example 1} (i).}\label{tab1}
{\scriptsize
{\renewcommand{\arraystretch}{1.4}
\begin{tabular}{crlllllllllll}
\toprule
\multirow{2}{*}{$\alpha$} & \multicolumn{1}{c}{\multirow{2}{*}{$\theta$}} & \multicolumn{5}{c}{Corrected scheme (\ref{novel.4})}                                                                                                                                  & \multicolumn{1}{c}{} & \multicolumn{5}{c}{Standard scheme (\ref{novel.3})}                                                                                                                                   \\ \cline{3-7} \cline{9-13}
                          & \multicolumn{1}{c}{}                          & \multicolumn{1}{c}{$\tau=2^{-5}$} & \multicolumn{1}{c}{$\tau=2^{-6}$} & \multicolumn{1}{c}{$\tau=2^{-7}$} & \multicolumn{1}{c}{$\tau=2^{-8}$} & \multicolumn{1}{c}{Rates} & \multicolumn{1}{c}{} & \multicolumn{1}{c}{$\tau=2^{-5}$} & \multicolumn{1}{c}{$\tau=2^{-6}$} & \multicolumn{1}{c}{$\tau=2^{-7}$} & \multicolumn{1}{c}{$\tau=2^{-8}$} & \multicolumn{1}{c}{Rates} \\ \midrule
\multirow{4}{*}{0.1}      & -0.9                                          & 4.33E-06                         & 3.10E-06                         & 6.92E-07                         & 1.62E-07                         & 2.09                      &                      & 7.50E-04                         & 3.91E-04                         & 1.96E-04                         & 9.82E-05                         & 1.00                      \\
                          & -0.5                                          & 1.86E-06                         & 8.76E-07                         & 2.65E-07                         & 7.13E-08                         & 1.89                      &                      & 1.86E-06                         & 8.76E-07                         & 2.65E-07                         & 7.13E-08                         & 1.89                      \\
                          & 0.5                                           & 1.47E-04                         & 3.43E-05                         & 8.27E-06                         & 2.02E-06                         & 2.03                      &                      & 2.02E-03                         & 9.97E-04                         & 4.95E-04                         & 2.47E-04                         & 1.01                      \\
                          & 0.9                                           & 2.53E-04                         & 5.78E-05                         & 1.38E-05                         & 3.37E-06                         & 2.03                      &                      & 2.87E-03                         & 1.41E-03                         & 6.95E-04                         & 3.46E-04                         & 1.01                      \\ \hline
\multirow{4}{*}{0.5}      & -0.8                                          & 1.15E-04                         & 2.49E-05                         & 5.78E-06                         & 1.39E-06                         & 2.05                      &                      & 3.15E-03                         & 1.60E-03                         & 8.04E-04                         & 4.03E-04                         & 1.00                      \\
                          & -0.5                                          & 3.86E-05                         & 6.97E-06                         & 1.44E-06                         & 3.24E-07                         & 2.15                      &                      & 3.86E-05                         & 6.97E-06                         & 1.44E-06                         & 3.24E-07                         & 2.15                      \\
                          & 0                                             & 2.35E-04                         & 5.70E-05                         & 1.40E-05                         & 3.49E-06                         & 2.01                      &                      & 5.49E-03                         & 2.72E-03                         & 1.35E-03                         & 6.74E-04                         & 1.00                      \\
                          & 0.6                                           & 2.35E-04                         & 5.70E-05                         & 1.40E-05                         & 3.49E-06                         & 2.01                      &                      & 1.23E-02                         & 6.02E-03                         & 2.98E-03                         & 1.49E-03                         & 1.01                      \\ \hline
\multirow{4}{*}{0.9}      & -0.5                                          & 2.35E-04                         & 5.70E-05                         & 1.40E-05                         & 3.49E-06                         & 2.01                      &                      & 3.05E-04                         & 7.23E-05                         & 1.76E-05                         & 4.35E-06                         & 2.02                      \\
                          & -0.2                                          & 1.28E-04                         & 2.95E-05                         & 7.10E-06                         & 1.74E-06                         & 2.03                      &                      & 6.78E-03                         & 3.30E-03                         & 1.63E-03                         & 8.10E-04                         & 1.01                      \\
                          & 0.3                                           & 3.56E-04                         & 8.65E-05                         & 2.14E-05                         & 5.31E-06                         & 2.01                      &                      & 1.78E-02                         & 8.72E-03                         & 4.33E-03                         & 2.15E-03                         & 1.01                      \\
                          & 0.6                                           & 7.64E-04                         & 1.84E-04                         & 4.51E-05                         & 1.12E-05                         & 2.01                      &                      & 2.44E-02                         & 1.20E-02                         & 5.95E-03                         & 2.96E-03                         & 1.01                      \\ \bottomrule
\end{tabular}}}
\end{table}
\begin{table}[]
\centering
\caption{$L^2$ error and convergence rates at time $t=0.5$ of \textit{Example 1} (ii).}\label{tab2}
{\scriptsize
{\renewcommand{\arraystretch}{1.4}
\begin{tabular}{crlllllllllll}
\toprule
\multirow{2}{*}{$\alpha$} & \multicolumn{1}{c}{\multirow{2}{*}{$\theta$}} & \multicolumn{5}{c}{Corrected scheme}                                                                                                                                      & \multicolumn{1}{c}{} & \multicolumn{5}{c}{Standard scheme}                                                                                                                                       \\ \cline{3-7} \cline{9-13}
                          & \multicolumn{1}{c}{}                          & \multicolumn{1}{c}{$\tau=2^{-5}$} & \multicolumn{1}{c}{$\tau=2^{-6}$} & \multicolumn{1}{c}{$\tau=2^{-7}$} & \multicolumn{1}{c}{$\tau=2^{-8}$} & \multicolumn{1}{c}{Rates} & \multicolumn{1}{c}{} & \multicolumn{1}{c}{$\tau=2^{-5}$} & \multicolumn{1}{c}{$\tau=2^{-6}$} & \multicolumn{1}{c}{$\tau=2^{-7}$} & \multicolumn{1}{c}{$\tau=2^{-8}$} & \multicolumn{1}{c}{Rates} \\ \midrule
\multirow{4}{*}{0.2}      & -0.5                                          & 2.68E-06                          & 7.74E-07                          & 2.03E-07                          & 5.14E-08                          & 1.98                      &                      & 2.68E-06                          & 7.74E-07                          & 2.03E-07                          & 5.14E-08                          & 1.98                      \\
                          & -0.3                                          & 7.66E-06                          & 1.92E-06                          & 4.80E-07                          & 1.18E-07                          & 2.02                      &                      & 9.41E-05                          & 4.69E-05                          & 2.28E-05                          & 1.07E-05                          & 1.09                      \\
                          & 0                                             & 1.83E-05                          & 4.39E-06                          & 1.07E-06                          & 2.62E-07                          & 2.03                      &                      & 2.42E-04                          & 1.19E-04                          & 5.75E-05                          & 2.68E-05                          & 1.10                      \\
                          & 0.9                                           & 7.69E-05                          & 1.75E-05                          & 4.14E-06                          & 9.97E-07                          & 2.06                      &                      & 7.07E-04                          & 3.40E-04                          & 1.63E-04                          & 7.56E-05                          & 1.11                      \\ \hline
\multirow{4}{*}{0.8}      & -0.5                                          & 8.79E-05                          & 2.12E-05                          & 5.20E-06                          & 1.28E-06                          & 2.03                      &                      & 8.79E-05                          & 2.12E-05                          & 5.20E-06                          & 1.28E-06                          & 2.03                      \\
                          & 0.1                                           & 1.99E-04                          & 4.64E-05                          & 1.12E-05                          & 2.71E-06                          & 2.04                      &                      & 7.59E-04                          & 3.95E-04                          & 1.95E-04                          & 9.18E-05                          & 1.09                      \\
                          & 0.5                                           & 3.28E-04                          & 7.47E-05                          & 1.77E-05                          & 4.27E-06                          & 2.05                      &                      & 1.36E-03                          & 6.82E-04                          & 3.31E-04                          & 1.54E-04                          & 1.10                      \\
                          & 0.7                                           & 4.11E-04                          & 9.26E-05                          & 2.18E-05                          & 5.25E-06                          & 2.06                      &                      & 1.68E-03                          & 8.29E-04                          & 3.99E-04                          & 1.86E-04                          & 1.10                      \\ \bottomrule
\end{tabular}}}
\end{table}
\\
\textit{Example 2.} We illustrate the robustness of (\ref{novel.4}) when $\alpha \to 0$.
Let $\Omega=(0,\pi), T=1$ and $u(x,t)=(E_\alpha(-t^\alpha)+t^3)\sin x$ such that $v=\sin x\in D(\Delta)$.
The source term is $f(x,t)=\big(6t^{3-\alpha}/\Gamma(4-\alpha)+t^3\big)\sin x$.
In Fig.\ref{C1} (a), we illustrate the $L^2$ error of the scheme (\ref{novel.4}) for varying $\alpha$ under different $\theta=-0.5,0.1,0.4,0.8$.
Particularly, the cases $\theta=0.1$ and $0.4$ of the scheme in \cite{YinLiuLiZhang1} are also presented.
Obviously, the scheme (\ref{novel.4}) is much more robust when $\alpha \to 0$ than the scheme in \cite{YinLiuLiZhang1}.
\par
It seems weird that in (\ref{novel.1.2}) the term $\phi(t_{n-\theta})$ is approximated by a nonlocal formula with coefficients $\theta_j$ with $j=0,1,\cdots, n$.
We shall argue that $\theta_j$ decays exponentially as plotted in Fig.\ref{C1} (b), and thus we only need the first few $\theta_j$'s.
\section{Conclusion}\label{sec.conc}
A general conversion strategy is proposed to develop robust and accurate difference formulas based on known ones by involving a shifted parameter $\theta$.
As a demonstration, the well-known BDF2 is considered and is proved rigorously for the subdiffusion problem (\ref{I.1}) showing that our scheme is robust even for very small $\alpha$ and can resolve the initial singularity of the solution.
\section*{Acknowledgments}
The work of the second author was supported by the NSF of Inner Mongolia 2021BS01003,
the third author was supported in part by Grants NSFC 12061053 and the NSF of Inner Mongolia 2020MS01003,
and the fourth author was supported in part by the grant NSFC 12161063 and the NSF of Inner Mongolia 2021MS01018.


\end{document}